\documentclass[12pt,dvips]{article}
\usepackage{amsmath,amsfonts,amssymb,amsthm,graphicx,verbatim,subfig}
\usepackage[all]{xy}

\ifx\pdftexversion\undefined

\usepackage[a4paper,colorlinks,link
=black,filecolor=black,citecolor=black,urlcolor=black,pdfstartview=FitH]{hyperref}
\else

\usepackage[a4paper,colorlinks,linkcolor=black,filecolor=black,citecolor=black,urlcolor=black,pdfstartview=FitH]{hyperref}
\fi

\font\sixbb=msbm6
\font\eightbb=msbm8
\font\twelvebb=msbm10 scaled 1095
\newfam\bbfam
\textfont\bbfam=\twelvebb \scriptfont\bbfam=\eightbb
                           \scriptscriptfont\bbfam=\sixbb

\newcommand{\Rea}{{\mathbb R}}
\newcommand{\Com}{{\mathbb C}}
\newcommand{\Int}{{\mathbb Z}}
\newcommand{\Rat}{{\mathbb Q}}
\newcommand{\Nat}{{\mathbb N}}
\newcommand{\FF}{{\mathbb F}}

\newcommand{\PP}{{\mathbb P}}

\newtheorem{theorem}{\bf Theorem}[section]
\newtheorem{claim}[theorem]{\bf Claim}

\newtheorem{proposition}[theorem]{\bf Proposition}

\newcommand{\enp}{\begin{flushright} $\Box$ \end{flushright}}
\newcommand{\beq}[0]{\begin{equation}}
\newcommand{\enq}[0]{\end{equation}}

\newcommand{\delp}{\Delta_{p-1}}
\newcommand{\ca}{{\cal A}}
\newcommand{\cb}{{\cal B}}
\newcommand{\cgg}{{\cal G}}

\newcommand{\cs}{{\cal S}}
\newcommand{\cp}{{\cal P}}
\newcommand{\crr}{{\cal R}}
\newcommand{\chh}{{\cal H}}
\newcommand{\rk}{{\rm rank}}
\newcommand{\thh}{\tilde{H}}
\newcommand{\supp}{{\rm supp}}
\newcommand{\cf}{{\cal F}}
\newcommand{\sgn}{{\rm sgn}}

\newcommand{\namedref}[2]{\hyperref[#2]{#1~\ref*{#2}}}

\title{Sum Complexes and Uncertainty Numbers}
\begin{document}
\author{Roy Meshulam
\thanks{Department of Mathematics,
Technion, Haifa 32000, Israel. e-mail:
meshulam@math.technion.ac.il~.}}
\maketitle
\pagestyle{plain}
\begin{abstract}
Let $p$ be a prime and let $A \subset \FF_p$. For $k<p$ let
$X_{A,k}$ be the $(k-1)$-dimensional complex on the vertex
set $\FF_p$ with a full $(k-2)$-skeleton whose $(k-1)$-faces are $\sigma \subset \FF_p$ such that
$|\sigma|=k$ and $\sum_{x \in \sigma}x \in A$. \\
The homology groups of $X_{A,k}$ with field coefficients are determined. In particular it is shown that
if $|A| \leq k$ then $H_{k-1}(X_{A,k};\FF_p)=0.$
\\
This implies a homological characterization
of uncertainty numbers of subsets of $\FF_p$: Let $\FF$ be algebraically closed, then
$H_{k-1}(X_{A,k};\FF)\neq 0$ iff there exists an $0 \neq f \in \FF[\FF_p]$ with $\supp(f) \subset A$ such that the endomorphism
of $\FF[\FF_p]$ given by $g \rightarrow fg$ has rank at most $p-k$.
\end{abstract}

\section{Introduction}
\label{intro}
\ \ \
Let $p$ be a prime and let $\FF_p=\{0,\ldots,p-1\}$ be the field of order $p$. Denote by $\delp$ the $(p-1)$-simplex
on the vertex set $\FF_p$ and let $\delp^{(k)}$ be the $k$-dimensional skeleton of $\delp$.
Let $A=\{a_1,\ldots,a_m\}$ be a subset of  $\FF_p$ and let $1<k<p$.
The {\it Sum Complex} $X_{A,k}$ was defined in \cite{LMR10} by
$$X_{A,k}=\delp^{(k-2)}~ \cup ~\{ \sigma \subset \FF_p~:~ |\sigma|=k~,~\sum_{x \in \sigma} x \in A\}.$$
$X_{A,k}$ is a $(k-1)$-dimensional complex whose $f$-vector satisfies
$f_i(X_{A,k})=\binom{p}{i+1}$ for $0 \leq i \leq k-2$ and $f_{k-1}(X_{A,k})=\frac{m}{p}\binom{p}{k}=\frac{m}{k}\binom{p-1}{k-1}$.
\ \\ \\
{\bf Example:} Let $p=7$ and $A=\{0,1,3\} \subset
\FF_7$. Then $X_{A,3}$ is homotopy equivalent to the real projective plane $\Rea\PP^2$,
see Figure 1 in \cite{LMR10}.




In this paper we study the homology of $X_{A,k}$ with field coefficients. Let $\FF$ be a field of characteristic $\ell$. First suppose that $\ell \nmid p$ and let $\omega$ be a primitive $p$-th root of unity in the algebraic closure $\overline{\FF}$. For $\beta=(b_1,\ldots,b_k) \in \FF_p^k$ let $M_{\beta}$ be the $k \times m$ matrix given by $M_{\beta}(i,j)=\omega^{b_i a_j}$. Let $$\cb_k=\{(b_1,\ldots,b_k):0 \leq b_1 < \cdots < b_k \leq p-1\}.$$
The case $m=k$ of the following result is implicit in \cite{LMR10}.
\begin{theorem}
\label{coprime}
If ${\rm char} \, \FF \nmid p$ then
$$
\dim H_{k-1}(X_{A,k};\FF)=\frac{m}{k}\binom{p-1}{k-1}-\frac{1}{p} \sum_{\beta \in \cb_k} \rk \, M_{\beta}.
$$
\end{theorem}
Our main result concerns the homology of $X_{A,k}$ with $\FF_p$ coefficients.
\begin{theorem}
\label{ptree}
\begin{equation}
\label{hmodp}
\dim H_{k-1}(X_{A,k};\FF_p)=
\left\{
\begin{array}{cl}
0 & m \leq k \\
(\frac{m}{k}-1)\binom{p-1}{k-1} & m > k.
\end{array}
\right.~~
\end{equation}
\end{theorem}
\noindent
{\bf Remarks:}
\\
1) The argument given in \cite{LMR10} for the case $m=k$ of Theorem \ref{coprime} does not seem to extend to the modular case. The approach here is different and is also used in the proof of our main result Theorem \ref{ptree}.
\\
2) The reduced Euler characteristic of $X_{A,k}$ is
\begin{equation}
\label{eulerc}
\begin{split}
\tilde{\chi}(X_{A,k}) &= -1+\sum_{i=0}^{k-2}(-1)^i \binom{p}{i+1}+(-1)^{k-1}\frac{m}{p}\binom{p}{k} \\
&= (-1)^{k-1}(\frac{m}{k}-1)\binom{p-1}{k-1}.
\end{split}
\end{equation}
Since $\thh_i(X_{A,k};\FF_p)=0$ for $0 \leq i <k-2$ it follows that
\begin{equation}
\label{hmodpa}
\begin{split}
\dim \thh_{k-2}(X_{A,k};\FF_p)&=\dim H_{k-1}(X_{A,k};\FF_p)-(\frac{m}{k}-1)\binom{p-1}{k-1} \\
&= \left\{
\begin{array}{cl}
(1-\frac{m}{k})\binom{p-1}{k-1} & m \leq k \\
0 & m > k
\end{array}
\right.~~
\end{split}
\end{equation}
\\
3) A classical result of Chebotar\"{e}v (see e.g. \cite{chebotarev})
asserts that for $\FF=\Rat$ all $M_{\beta}$'s have full rank. Theorem \ref{coprime} therefore implies
that (\ref{hmodp}) and (\ref{hmodpa}) remain true for $\thh_*(X_{A,k};\Rat)$.

We next discuss an application of Theorems \ref{coprime} and \ref{ptree} to discrete uncertainty principles.
Let $K$ be a finite  abelian group and let $\FF[K]$ be the group algebra of $K$ over the field $\FF$. For an element $f \in \FF[K]$ let  $T_f:\FF[K] \rightarrow \FF[K]$ be given by $T_f g=fg$.
The {\it uncertainty number} of a subset $A \subset K$ is defined by
$$u_{\FF}(A)=\min\{ \rk \, T_f~:~ \emptyset \neq \supp(f)  \subset A\}.$$

The motivation for this definition and terminology is as follows. Let $m$ be the exponent of $K$ and
suppose $\FF$ contains a primitive $m$-th root of unity. Let $\widehat{K}$ denote the group of $\FF$-valued characters of $K$. Identifying $\FF[K]$ with the space of $\FF$-valued functions on $K$, the Fourier Transform of a function $f \in \FF[K]$ is the function $\widehat{f} \in \FF[\widehat{K}]$ given by $\widehat{f}(\chi)=\sum_{x \in K} \chi(-x) f(x)$.
It is well known that $\rk \, T_f=|\supp(\widehat{f})|$. Thus in the semisimple case
\begin{equation}
\label{unfour}
u_{\FF}(A)=\min\{|\supp(\widehat{f})|~:~ \emptyset \neq \supp(f)  \subset A\}.
\end{equation}

The classical discrete uncertainty principle (see e.g. \cite{DS89}) asserts that  $u_{\FF}(A) \geq \frac{|K|}{|A|}$
for any $\FF$ and $\emptyset \neq A \subset K$.  While this bound is sharp when $A$ is a coset of $K$, it can often be improved for particular choices of $K,A$ and $\FF$. One such example is a
result of Tao \cite{Tao05} asserting that if $A \subset \FF_p$ then $u_{\Com}(A)=p-|A|+1$. See \cite{M06} for an extension to general abelian groups. Here we note a simple relation between the homology of sum complexes and uncertainty numbers of subsets of $\FF_p$.
\begin{theorem}
\label{unhom}
If $\FF$ is algebraically closed then for any $A \subset \FF_p$
$$u_{\FF}(A)=p-\max\{k:H_{k-1}(X_{A,k};\FF) \neq 0\}.$$
\end{theorem}
\noindent
{\bf Example:} It is easy to see that $u_{\FF}(A) \geq p- \max A$ for any
$A \subset \FF_p$ and any field $\FF$. Taking $A=\{0,1,3\} \subset \FF_7$ it follows
that $u_{\overline{\FF_2}}(A) \geq 4$. On the other hand,
as noted earlier $X_{A,3}$ is homotopy equivalent to $\Rea\PP^2$, hence $H_2(X_{A,3};\FF_2) \neq 0$. Theorem \ref{unhom} then implies that $u_{\overline{\FF_2}}(A)=4$.
It can be checked that in fact
$u_{\FF_2}(A)=4$.
\ \\

Let $C_p$ be the multiplicative cyclic group of order $p$ and let $G=C_p^k$.
In Section \ref{s:cycles} we identify $H_{k-1}(X_{A,k};\FF)$ with a certain subspace $\chh(A)$ of skew-symmetric
elements of the group algebra $\FF[G]$. This characterization is used in Section \ref{s:semi}
to prove Theorem  \ref{coprime}. The proof of Theorem \ref{ptree} given in Section \ref{s:modu} is more involved and depends additionally on some properties of generalized Vandermonde determinants over the group algebra $\FF_p[G]$.  Theorem \ref{unhom} is derived in Section \ref{s:unh} as a  direct consequence of Theorems \ref{coprime} and \ref{ptree}. We conclude in Section \ref{s:remarks} with some comments and open problems.

\section{A Characterization of Cycles}
\label{s:cycles}
\ \ \
Let $\cf(\FF_p^k,\FF)$ denote the space of $\FF$-valued functions on $\FF_p^{k}$.
A function $\phi \in \cf(\FF_p^k,\FF)$ is {\it skew-symmetric} if
$\phi(\gamma_{\sigma^{-1}(1)},\ldots,\gamma_{\sigma^{-1}(k)})=\sgn(\sigma)\phi(\gamma_1,\ldots,\gamma_k)$
for all $(\gamma_1,\ldots,\gamma_k) \in \FF_p^k$ and $\sigma$ in the symmetric group $S_k$. If ${\rm char} \,\FF=2$ then $\phi$ is additionally required to satisfy
 $\phi(\gamma_1,\ldots,\gamma_k)=0$
if $\gamma_i=\gamma_j$ for some $i \neq j$.
Let $\ca(\FF_p^{k},\FF)$ be the space of skew-symmetric functions in $\cf(\FF_p^k,\FF)$.

Fix a generating set $\{x_1,\ldots,x_k\}$ of $G$ and let $x=(x_1,\ldots,x_k) \in G^k$.
For $\gamma=(\gamma_1,\ldots,\gamma_k) \in \FF_p^k$
we abbreviate $x^{\gamma}=\prod_{j=1}^k x_j^{\gamma_j}$. Let $q:\cf(\FF_p^k,\FF) \rightarrow \FF[G]$ be the $\FF$-linear isomorphism given by
$$q(\phi)=\sum_{\gamma \in \FF_p^k} \phi(\gamma) x^{\gamma}.$$
An element $s \in \FF[G]$ is {\it homogenous of degree} $d \in \FF_p$ if $s=q(\phi)$
where $$\supp(\phi) \subset W_d=\{\alpha=(\alpha_1,\ldots,\alpha_k) \in \FF_p^k: \sum_{i=1}^k \alpha_i=d\}.$$
Let $\FF[G]_d$ denote the space of degree $d$ elements of $\FF[G]$
and let $\rho_d$ be the projection from $\FF[G]$ onto $\FF[G]_d$.
An element $s \in \FF[G]$ is {\it skew-symmetric} if $s=q(\phi)$ for some $\phi \in \ca(\FF_p^{k},\FF)$. Let $\cs$ be the space of skew-symmetric elements of $\FF[G]$ and let $\cs_d=\cs \cap \FF[G]_d$.

The space of $\FF$-valued $(k-1)$-chains of $X_{A,k}$ is
$$C_{k-1}(X_{A,k};\FF)=\{\phi \in \ca(\FF_p^{k},\FF)~:~\supp(\phi) \subset \cup_{a \in A} W_a\}.$$
A $(k-1)$-chain $\phi \in C_{k-1}(X_{A,k};\FF)$ is a $(k-1)$-cycle of $X_{A,k}$ if
\begin{equation}
\label{kcycle}
\sum_{a \in A} \phi(\alpha_1,\ldots,\alpha_{i-1},a-\sum_{\{j:j \neq i\}} \alpha_j,\alpha_{i+1},\ldots, \alpha_k)=0
\end{equation}
for all fixed $1 \leq i \leq k$
and $(\alpha_1,\ldots,\widehat{\alpha_i},\ldots,\alpha_k) \in \FF_p^{k-1}$.
Let
$$
\chh(A)=\{s \in \bigoplus_{a \in A} \cs_a~:~\sum_{a \in A}x_i^{-a} \rho_a(s)=0 {\rm ~~for~all~~} 1 \leq i \leq k\}.$$
The homology space $H_{k-1}(X_{A,k};\FF)=Z_{k-1}(X_{A,k};\FF)$ is characterized by the following
\begin{claim}
\label{mainid}
$$q \left( H_{k-1}(X_{A,k};\FF) \right)=\chh(A).$$
\end{claim}
\noindent
{\bf Proof:} Let $\phi \in C_{k-1}(X_{A,k};\FF)$ and fix an $1 \leq i \leq k$. Then
\begin{equation*}
\begin{split}
& \sum_{a \in A} x_i^{-a} \rho_a(q(\phi))=\sum_{a \in A}  x_i^{-a} \sum_{(\alpha_1,\ldots,\alpha_k) \in W_a}
\phi(\alpha_1,\ldots,\alpha_k) \prod_{j=1}^k x_j^{\alpha_j} \\
& = \sum_{a \in A}  x_i^{-a} \sum_{(\alpha_j:j \neq i) \in \FF_p^{k-1}}
\phi(\alpha_1,\ldots,\alpha_{i-1},a-\sum_{j \neq i} \alpha_j,\alpha_{i+1},\ldots, \alpha_k)
(\prod_{j \neq i}x_j^{\alpha_j}) x_i^{a-\sum_{j \neq i} \alpha_j} \\
& =
\sum_{(\alpha_j:j \neq i) \in \FF_p^{k-1}}
\left(
\sum_{a \in A} \phi(\alpha_1,\ldots,\alpha_{i-1},a-\sum_{j \neq i} \alpha_j,\alpha_{i+1},\ldots, \alpha_k)
\right)
\prod_{j \neq i} (x_jx_i^{-1})^{\alpha_j} =0.
\end{split}
\end{equation*}
Therefore $\sum_{a \in A} x_i^{-a} \rho_a(q(\phi))=0$ iff for all $(\alpha_j:j \neq i) \in \FF_p^{k-1}$
$$\sum_{a \in A} \phi(\alpha_1,\ldots,\alpha_{i-1},a-\sum_{j \neq i} \alpha_j,\alpha_{i+1},\ldots, \alpha_k)=0.$$
Hence the Claim follows from (\ref{kcycle}).
{\enp}

\section{The Semisimple Case}
\label{s:semi}

In this section we prove Theorem \ref{coprime}. We may assume that $\FF$ is algebraically closed.
Let
$$\crr=\{(r_a)_{a \in A} \in \cs^A~:~
\sum_{a \in A}x_i^{-a}r_a=0 {\rm ~~for~all~~} 1 \leq i \leq k\}.$$
Let $y=\prod_{j=1}^k x_j$ and let $P \subset G$ be the cyclic group generated by $y$.
\begin{claim}
\label{allsol}
The mapping $B:\FF[P] \otimes_{\FF} \chh(A) \rightarrow \crr$ given by
$$B(u \otimes s)=(u \rho_a(s))_{a \in A}$$
is an isomorphism.
\end{claim}
\noindent
{\bf Proof:} We first show that $B$ is injective. Let $w=\sum_{j=0}^{p-1} y^j \otimes s_j \in \ker B$.
Then $\sum_{j=0}^{p-1} y^j \rho_a(s_j)=0$ for all $a \in A$.
Since $y^j \rho_a(s_j) \in \cs_{a+kj}$ it follows that $y^j \rho_a(s_j)=0$ and hence $\rho_a(s_j)=0$ for all $0 \leq j \leq p-1$ and $a \in A$. Therefore $w=0$.
To show surjectivity let $(r_a)_{a \in A} \in \crr$.
For $0 \leq i,j \leq p-1$ let
$$s_j=y^{-j} \sum_{a \in A} \rho_{a+kj}(r_a) \in \bigoplus_{a \in A} \cs_a$$
and
$$t_{ij} =\sum_{a \in A} x_i^{-a} \rho_{a+kj}(r_a) \in \cs_{kj}.$$
Then for any $0 \leq i \leq p-1$
$$
0=\sum_{a \in A} x_i^{-a}r_a=\sum_{a \in A}x_i^{-a} \sum_{j=0}^{p-1} \rho_{a+kj}(r_a)
=\sum_{j=0}^{p-1} t_{ij}.
$$
It follows that $t_{ij}=0$ for all $0 \leq i,j \leq p-1$. Therefore
$$
\sum_{a \in A} x_i^{-a} \rho_a(s_j)=\sum_{a \in A} x_i^{-a} y^{-j} \rho_{a+kj}(r_a)=y^{-j} t_{ij}=0$$
and hence $s_j \in \chh(A)$. Finally
$$\sum_{j=0}^{p-1} y^j \rho_a(s_j)=\sum_{j=0}^{p-1} y^j (y^{-j}\rho_{a+kj}(r_a))=
\sum_{j=0}^{p-1}\rho_{a+jk}(r_a)=r_a ~,$$
therefore $B(\sum_{j=0}^{p-1} y^j \otimes s_j)=(r_a)_{a \in A}.$
{\enp}
\noindent
Claim \ref{allsol} implies that $\dim \crr=p \, \dim \chh(A).$
Theorem \ref{coprime} will thus follow from
\begin{proposition}
\label{dimsall}
$$
\dim \crr=m\binom{p}{k}-\sum_{\beta \in \cb_k} \rk \, M_{\beta}.
$$
\end{proposition}
\noindent
{\bf Proof:}  Recall that $\omega$ is a primitive $p$-th root of unity in $\FF=\overline{\FF}$. The Fourier Transform is the automorphism of $\cf(\FF_p^k,\FF)$
given by
$$\widehat{\phi}(\beta)=\sum_{\alpha \in \FF_p^k} \phi(\alpha) \omega^{-\beta \alpha}.$$
Define an $\FF$-linear isomorphism
$$\Phi: \cs^A \rightarrow \bigoplus_{\beta \in \cb_k} \FF^m$$
as follows. For $r=(r_a)_{a \in A} \in \cs^A$ where $r_a=q(\phi_a)$
and $\phi_a \in \ca(\FF_p^k,\FF)$ let
$$\Phi(r)=\left( (\widehat{\phi_{a_j}}(\beta))_{j=1}^m : \beta  \in \cb_k \right).$$
\begin{claim}
\label{isom}
$\Phi$ restricts to an isomorphism from $\crr$ onto $\oplus_{\beta \in \cb_k} \ker M_{\beta}$.
\end{claim}
\noindent
{\bf Proof:} For $1 \leq i \leq k$ let $e_i$ be the $i$-th unit vector in $\FF_p^k$. Then
$\psi_i=q^{-1}(\sum_{a \in A}x_i^{-a}r_a) \in \cf(\FF_p^k,\FF)$ is given by
$$\psi_i(\alpha)=\sum_{a \in A} \phi_a(\alpha+a e_i).$$
For $\beta=(b_1,\ldots,b_k) \in \FF_p^k$
\begin{equation*}
\begin{split}
\widehat{\psi_i}(\beta)&=\sum_{\alpha \in \FF_p^k} \sum_{a \in A} \phi_a(\alpha+ae_i) \omega^{-\beta \alpha}
= \sum_{a \in A} \sum_{\alpha \in \FF_p^k} \phi_a(\alpha) \omega^{-\beta(\alpha-ae_i)} \\
&=\sum_{a \in A} \omega^{b_i a} \sum_{\alpha \in \FF_p^k} \phi_a(\alpha) \omega^{-\beta\alpha}=
\sum_{a \in A} \omega^{b_i a} \widehat{\phi_a}(\beta).
\end{split}
\end{equation*}
It follows that $r=(r_a)_{a \in A} \in \crr$ iff $\psi_1=\cdots=\psi_k=0$ iff
$$
M_{\beta}\left[
\begin{array}{c}
\widehat{\phi_{a_1}}(\beta) \\
\vdots \\
\widehat{\phi_{a_m}}(\beta)
\end{array}
\right]=
\left[
\begin{array}{ccc}
\omega^{b_1 a_1} &  \cdots & \omega^{b_1 a_{m}} \\
\vdots & \ddots & \vdots \\
\omega^{b_k a_1} & \cdots & \omega^{b_k a_{m}}
\end{array}
\right]
\left[
\begin{array}{c}
\widehat{\phi_{a_1}}(\beta) \\
\vdots \\
\widehat{\phi_{a_m}}(\beta)
\end{array}
\right]=
\left[
\begin{array}{c}
\widehat{\psi_{1}}(\beta) \\
\vdots \\
\widehat{\psi_{m}}(\beta)
\end{array}
\right]=0
$$
for all $\beta=(b_1,\ldots,b_k) \in \FF_p^k$. Therefore $\Phi(\crr) \subset \oplus_{\beta \in \cb_k} \ker M_{\beta}$.
The bijectivity of the restriction $\Phi_{|\crr}$ follows from
the bijectivity $\Phi$.
{\enp}
\noindent
By Claim \ref{isom}
\begin{equation*}
\begin{split}
\dim \crr &=\sum_{\beta \in \cb_k} \dim \ker M_{\beta} = \sum_{\beta \in \cb_k} (m-\rk \, M_{\beta}) \\
&= m\binom{p}{k}-\sum_{\beta \in \cb_k} \rk \, M_{\beta}.
\end{split}
\end{equation*}
{\enp}

\section{The Modular Case}
\label{s:modu}
\ \ \ \
In subsections  \ref{groupa} and
\ref{groupb} we study certain properties of determinants of generalized Vandermonde matrices over the group algebra $\FF_p[G]$. These results are then used in Section \ref{modph} to prove Theorem \ref{ptree}.

\subsection{A Generalized Vandermonde}
\label{groupa}
\ \ \ \
Recall that $\{x_1,\ldots,x_k\}$ is a fixed generating set of $G$ and $x=(x_1,\ldots,x_k)$.
For $\beta=(b_1,\ldots,b_k) \in \cb_k$ let
$$
N_{\beta}=
\left[
\begin{array}{ccc}
x_1^{-b_1} & \cdots & x_1^{-b_{k}} \\
\vdots &  \ddots & \vdots \\
x_k^{-b_1} & \cdots & x_k^{-b_{k}}
\end{array}
\right].
$$
\begin{proposition}
\label{schur}
\begin{equation}
\label{detschur}
\det N_{\beta} =w_{\beta} \prod_{1 \leq i< j \leq k} (x_i-x_j)
\end{equation}
where $w_{\beta}$ is a unit of $\FF_p[G]$.
\end{proposition}

Recall the definition of Schur polynomials (see e.g. \cite{E11}).
Let $\xi=(\xi_1,\ldots,\xi_k)$ be a vector of variables.
For a partition $\lambda=(\lambda_1 \geq \cdots \geq \lambda_k)$ let
$$D_\lambda(\xi)=D_{\lambda}(\xi_1,\ldots,\xi_k)=\det([\xi_i^{\lambda_j+k-j}]_{i,j=1}^k) \in \Int[\xi_1,\ldots,\xi_k].$$
Note that for the zero partition $0=(0,\ldots,0)$
$$D_0(\xi)=\det  \left[
\begin{array}{cccc}
\xi_1^{k-1} & \xi_1^{k-2} &  \cdots & 1 \\
\vdots &  \vdots & \ddots & \vdots \\
\xi_k^{k-1} & \xi_k^{k-2} &  \cdots & 1
\end{array}
\right]=
\prod_{1 \leq i<j \leq k} (\xi_i-\xi_j).$$
The Schur polynomial associated with $\lambda$ is
$$s_{\lambda}(\xi)=\frac{D_{\lambda}(\xi)}{D_0(\xi)} \in \Int[\xi_1,\ldots,\xi_k].$$
The dimension formula (see e.g. Proposition 5.21.2 in \cite{E11}) asserts that
\begin{equation}
\label{dimrep}
s_{\lambda}(1,\ldots,1)=\prod_{1 \leq i < j \leq k}\frac{\lambda_i-\lambda_j+j-i}{j-i}.
\end{equation}
\ \\ \\
{\bf Proof of Proposition \ref{schur}:} Let
$$\lambda=(\lambda_1,\ldots,\lambda_k)=(p-b_1-k+1,p-b_2-k+2,\ldots,p-b_k).$$
Then
\begin{equation*}
\begin{split}
\det N_{\beta} &= D_{\lambda}(x)=s_{\lambda}(x) D_0(x) \\
&= s_{\lambda}(x) \prod_{1 \leq i< j \leq k} (x_i-x_j).
\end{split}
\end{equation*}
By (\ref{dimrep}) the image of $s_{\lambda}(x)\in \FF_p[G]$ under the augmentation map $\FF_p[G] \rightarrow \FF_p$ is
$$s_{\lambda}(1,\ldots,1)\, ({\rm mod} \, p)=\prod_{1 \leq i<j \leq k} \frac{b_j-b_i}{j-i}\, ({\rm mod} \, p) \neq ~0\, ({\rm mod} \, p).$$
It follows that $w_{\beta}=s_{\lambda}(x)$ is invertible in $\FF_p[G]$.
{\enp}

\subsection{Skew-Symmetric Annihilators of $D_0(x)$}
\label{groupb}
\ \ \ \
In this subsection we show
\begin{proposition}
\label{qskew}
Let $s \in \cs$. If $D_0(x)s=0$ then $s=0$.
\end{proposition}

The proof of Proposition \ref{qskew} depends on Proposition \ref{combs} below.
Let $\Nat$ denote the nonnegative integers and for $a,b \in \Nat$ let $[a,b]=\{a,\ldots,b\}$.
Let
$$\Nat_k=\{(\alpha_1,\ldots,\alpha_k) \in \Nat^k~:~ \alpha_i \neq \alpha_j ~~{\rm if}~~ i \neq j\}.$$
For
$\alpha=(\alpha_1,\ldots,\alpha_k)~,~ \beta=(\beta_1,\ldots,\beta_k) \in \Nat_k$ write
$\alpha \preceq \beta$ if
$\{\alpha_1,\ldots,\alpha_k\}$ precedes  $\{\beta_1,\ldots,\beta_k\}$
in the lexicographic order on subsets of $\Nat$, i.e. if
$$\sum_{i=1}^k 2^{-\alpha_i} \geq \sum_{i=1}^k 2^{-\beta_i}.$$
Fix an $\alpha=(\alpha_1,\ldots,\alpha_k) \in \Nat_k$ such that $\alpha_1<\cdots <\alpha_k$ and
let
$$L=\{1 \leq i \leq k-1: \alpha_i+1<\alpha_{i+1}\}.$$
Write $L=\{\ell_1<\ldots<\ell_{t-1}\}$ and let $\ell_0=0~,~ \ell_t=k$.
For $1 \leq i \leq t$ let $A_i=[\ell_{i-1}+1,\ell_i]$.
Let
$$\cgg_1(\alpha)=\{(\gamma,\sigma)\in \Nat_k \times S_k~:~\gamma \preceq \alpha~{\rm and}~ \gamma_j-\sigma(j)=\alpha_j-j ~~{\rm for~all}~~ j\}.$$
We'll need the following characterization of $\cgg_1(\alpha)$.
Let $S_A$ denote the symmetric group on a set $A$. Let $T=S_{A_1} \times \cdots \times S_{A_t}$ be the Young subgroup of $S_k$ corresponding to the partition $[k]= \cup_{i=1}^t A_i$. Let
$$\cgg_2(\alpha)=\{(\gamma,\sigma)\in \Nat_k \times T~:~
\gamma_j=\alpha_{\sigma(j)}~~{\rm for~all}~~j\}.$$
\begin{proposition}
\label{combs}
$\cgg_1(\alpha)=\cgg_2(\alpha)$
\end{proposition}
\noindent
{\bf Proof:}  We first show that $\cgg_2(\alpha) \subset \cgg_1(\alpha)$. Let $(\gamma,\sigma) \in \cgg_2(\alpha)$ and
let $1 \leq j \leq k$. If
$j \in A_i$ then $\sigma(j) \in A_i$ and  hence $\alpha_{\sigma(j)}-\alpha_j=\sigma(j)-j$. Therefore
\begin{equation*}
\begin{split}
\gamma_j-\sigma(j)= \alpha_{\sigma(j)}-\sigma(j)=\alpha_j-j
\end{split}
\end{equation*}
and so $(\gamma,\sigma) \in \cgg_1(\alpha)$.
For the other direction
let $(\gamma,\sigma)\in \cgg_1(\alpha)$. Write $\gamma=(\gamma_1,\ldots,\gamma_k)$ and let $\pi \in S_k$ such that
$\gamma_{\pi(1)}< \cdots <\gamma_{\pi(k)}$.
\begin{claim}
\label{firsta}
For $1 \leq i \leq t$ and $j \in A_i$
\\
(a) $\sigma(\pi(j))=j$.
\\
(b) $\gamma_{\pi(j)}=\alpha_j$.
\\
(c) $\pi(j) \in A_i$.
\end{claim}
\noindent
{\bf Proof:} We argue by induction on $j$. Suppose (a),(b) and (c) hold for all $j'<j$.
(a) implies that $\{\sigma(\pi(j')): j'<j\}=[j-1]$ and hence $\sigma(\pi(j)) \geq j$.
Therefore
\begin{equation}
\label{one}
\alpha_j-j \geq \alpha_j -\sigma(\pi(j)).
\end{equation}
Next note that by (b) $\gamma_{\pi(j')}=\alpha_{j'}$ for all $j'<j$. As $\gamma \preceq \alpha$
it follows that $\alpha_j \geq \gamma_{\pi(j)}$ and therefore
\begin{equation}
\label{two}
\alpha_j -\sigma(\pi(j)) \geq \gamma_{\pi(j)}-\sigma(\pi(j))=\alpha_{\pi(j)}-\pi(j).
\end{equation}
Finally (c) implies that $\{\pi(j'):1 \leq j' \leq \ell_{i-1}\}=[1,\ell_{i-1}]$ and therefore
$\pi(j) \geq \ell_{i-1}+1$. Together with the assumption $j \in A_i$ it follows that
\begin{equation}
\label{three}
\alpha_{\pi(j)}-\pi(j) \geq \alpha_{\ell_{i-1}+1}-(\ell_{i-1}+1)=\alpha_j-j.
\end{equation}
It follows that the three inequalities in
(\ref{one}),(\ref{two}),(\ref{three}) are in fact equalities.
Therefore $\sigma(\pi(j))=j$, $\gamma_{\pi(j)}=\alpha_j$ and
$\alpha_{\pi(j)}=\alpha_j+(\pi(j)-j)$ respectively establishing (a),(b),(c) for $j$.
{\enp}
Claim 4.4 implies that $\sigma=\pi^{-1} \in T$ and that  $\gamma_{j}=\alpha_{\sigma(j)}$ for
all $1 \leq j \leq k$. Therefore $(\gamma,\sigma) \in  \cgg_2(\alpha)$.
{\enp}
\noindent
{\bf Proof of Proposition \ref{qskew}:} Let $s \in \cs$ such that $D_0(x)s=0$ and write $s=q(\phi)$ with
$\phi \in \ca(\FF_p^{k},\FF_p)$. We have to show that $\phi=0$.
By assumption
\begin{equation}
\label{rgamma}
\begin{split}
0 = D_0(x)s &=\sum_{\sigma\in S_k} \sgn(\sigma) \prod_{j=1}^k x_j^{k-\sigma(j)} \sum_{\gamma=(\gamma_1,\ldots,\gamma_k)\in \FF_p^k} \phi(\gamma) \prod_{j=1}^k x_j^{\gamma_j} \\
&=
\sum_{\gamma=(\gamma_1,\ldots,\gamma_k)\in \FF_p^k} \sum_{\sigma\in S_k} \sgn(\sigma)\phi(\gamma) \prod_{j=1}^k x_j^{\gamma_j+k-\sigma(j)}.
\end{split}
\end{equation}
Suppose for contradiction that $\phi \neq 0$ and let
$$\alpha=(\alpha_1,\ldots,\alpha_k)=\max\{\gamma \in \cb_k~:~\phi(\gamma) \neq 0\}$$
where the maximum is taken with respect to $\preceq$.
Let $\lambda \in \FF_p$ denote the coefficient of $\prod_{j=1}^k x_j^{\alpha_j+k-j}$ in
the expansion of $D_0(x)s~$ in the standard basis $\{x^{\beta}:\beta \in \FF_p^k\}$.
Note that if
$$\prod_{j=1}^k x_j^{\alpha_j+k-j}=\prod_{j=1}^k x_j^{\gamma_j+k-\sigma(j)}$$
then for all $1 \leq j \leq k$ $$\alpha_j-j = \gamma_j-\sigma(j) \, ({\rm mod} \, p).$$
Since
$$-1 \leq \alpha_j-j \leq p-1-k$$
and
$$-k \leq \gamma_j-\sigma(j) \leq p-2$$
it follows that $$\alpha_j-j = \gamma_j-\sigma(j).$$
Hence, Eq. (\ref{rgamma}) and Proposition \ref{combs} imply that
\begin{equation*}
\begin{split}
\lambda &=\sum_{(\gamma,\sigma) \in \cgg_1(\alpha)} \sgn(\sigma) \phi(\gamma)=
\sum_{(\gamma,\sigma) \in \cgg_2(\alpha)} \sgn(\sigma) \phi(\gamma) \\
&= \sum_{\sigma \in S_{A_1} \times \cdots \times S_{A_t}} \sgn(\sigma) \phi(\alpha_{\sigma(1)},\ldots,\alpha_{\sigma(k)}) \\
&= |S_{A_1} \times \cdots \times S_{A_t}| \, \phi(\alpha) = \prod_{i=1}^t (\ell_i-\ell_{i-1})! \, \phi(\alpha).
\end{split}
\end{equation*}
Since $\ell_t=k<p$ it follows that $\prod_{i=1}^t (\ell_i-\ell_{i-1})! \neq 0 ({\rm mod} \, p)$ and so $\lambda \neq 0$.  Therefore $D_0(x)s \neq 0$, a contradiction.
{\enp}

\subsection{Homology of $X_{A,k}$ over $\FF_p$}
\label{modph}

In this subsection we prove
Theorem \ref{ptree}. We first consider the case $m=k$.
\begin{theorem}
\label{ptreea}
If $|A|=k$ then
$H_{k-1}(X_{A,k};\FF_p)=0.$
\end{theorem}
\noindent
{\bf Proof:} Let $A=\{a_1,\ldots,a_k\}$ where $\alpha=(a_1,\ldots,a_k) \in \cb_k$.
Let $\phi \in H_{k-1}(X_{A,k};\FF_p)$ and let $s=q(\phi)$. By Claim \ref{mainid}
$$
N_{\alpha}
\left[
\begin{array}{c}
\rho_{a_1}(s) \\
\vdots \\
\rho_{a_{k}}(s)
\end{array}
\right]=
\left[
\begin{array}{ccc}
x_1^{-a_1} &  \cdots & x_1^{-a_{k}} \\
\vdots & \ddots & \vdots \\
x_k^{-a_1} & \cdots & x_k^{-a_{k}}
\end{array}
\right]
\left[
\begin{array}{c}
\rho_{a_1}(s) \\
\vdots \\
\rho_{a_{k}}(s)
\end{array}
\right]=0.
$$
Therefore $\, \det N_{\alpha} \cdot \rho_a(s)=0$ for all $a \in A$.
Proposition \ref{schur} implies that  $D_0(x)\rho_a(s)=0$
and hence $\rho_a(s)=0$ by Proposition \ref{qskew}. It follows that $\phi=0$ and so $H_{k-1}(X_{A,k};\FF_p)=0$.
{\enp}

\noindent
{\bf Proof of Theorem \ref{ptree}:} Let $|A|=m \geq k$ and let $A'$ be an arbitrary subset of $A$ of cardinality $k$.
Theorem \ref{ptree} and Eq. (\ref{eulerc}) imply that $\thh_*(X_{A',k};\FF_p)=0.$
Hence by the exact sequence
$$
0 =H_{k-1}(X_{A',k};\FF_p) \rightarrow H_{k-1}(X_{A,k};\FF_p) \rightarrow
$$
$$
\rightarrow H_{k-1}(X_{A,k},X_{A',k};\FF_p)
\rightarrow \thh_{k-2}(X_{A',k};\FF_p)=0
$$
it follows that
\begin{equation*}
\begin{split}
\dim H_{k-1}(X_{A,k};\FF_p) &=\dim H_{k-1}(X_{A,k},X_{A',k};\FF_p) \\
&= f_{k-1}(X_{A,k})-f_{k-1}(X_{A',k})=(\frac{m}{k}-1)\binom{p-1}{k-1}.
\end{split}
\end{equation*}
{\enp}

\section{Uncertainty Numbers and Homology}
\label{s:unh}
{\bf Proof of Theorem \ref{unhom}:} Let ${\rm char} \, \FF=\ell$ and consider two cases:
\ \\ \\
(i) $\ell \nmid p$. For $\lambda=(\lambda_1,\ldots,\lambda_m) \in \FF^m$ let
$f_{\lambda}:\FF_p \rightarrow \FF$ be given by $f_{\lambda}=\sum_{j=1}^m \lambda_j \delta_{a_j}$.
Then $\supp(f) \subset A$ and for $\beta=(b_1,\ldots,b_k) \in \FF_p^k$
\begin{equation}
\label{flambda}
M_{\beta} \lambda =(\sum_{j=1}^m \lambda_j \omega^{b_i a_j})_{i=1}^k= \left(\widehat{f_{\lambda}}(-b_i)\right)_{i=1}^k.
\end{equation}
Theorem \ref{coprime} implies that $H_{k-1}(X_{A,k};\FF) \neq 0$
iff $\rk \, M_{\beta}<m$ for some $\beta=(b_1,\ldots,b_k) \in \cb_k$.
On the other hand, (\ref{flambda}) implies that $\rk \, M_{\beta}<m$ iff there exists a $0 \neq \lambda \in \FF^m$
such $\widehat{f_{\lambda}}(-b_i)=0$ for $1 \leq i \leq k$.
It follows that $H_{k-1}(X_{A,k};\FF) \neq 0$ iff there exists a nonzero $f:\FF_p \rightarrow \FF$ such that
$\supp(f) \subset A$ and $|\supp(\widehat{f})| \leq p-k$.
\ \\ \\
(ii) $\ell=p$. Let $\FF$ be a field of characteristic $p$. By Theorem \ref{ptree}
$$p-\max\{k:H_{k-1}(X_{A,k};\FF) \neq 0\}=p-m+1.$$
It thus suffices to show that $u_{\FF}(A)=p-m+1$.
Let $f:\FF_p \rightarrow \FF$ be a function
with $\emptyset \neq \supp(f) \subset A=\{a_1,\ldots,a_m\}$ and let $g(x)=\sum_{i=1}^m f(a_i) x^{a_i} \in \FF[x]$.
Then
\begin{equation}
\label{rgcd}
\begin{split}
\rk \, T_f &=p-\deg \,{\rm gcd}(g(x),x^p-1) \\
&=p-\deg \, {\rm gcd}(g(x),(x-1)^p)=p-\mu(g)
\end{split}
\end{equation}
where $\mu(g)$ is the multiplicity of $1$ as a root of $g(x)$. By a result of Frenkel (Lemma 2 in \cite{Frenkel}),
$\mu(g) \leq m-1$ and hence $u_{\FF}(A) \geq p-m+1$. For the other direction note that the space of polynomials
$$\cp=\{g(x) \in \FF[x]:\deg \, g \leq p-1~,~\mu(g) \geq m-1\}$$
satisfies $\dim_{\FF} \cp = p-m+1$.  Hence  $\cp$ contains a nonzero $g$ of the form
$g(x)=\sum_{i=1}^m \lambda_i x^{a_i}$. It follows by (\ref{rgcd}) that $0 \neq f=\sum_{i=1}^m \lambda_i \delta_{a_i}:\FF_p \rightarrow \FF$
satisfies $\rk \,T_f \leq p-m+1$. Therefore $u_{\FF}(A) \leq p-m+1$.
{\enp}

\section{Concluding Remarks}
\label{s:remarks}
\ \ \
We mention two problems related to the results of this paper.
\begin{itemize}

\item
Let $X$ be a $(k-1)$-dimensional complex $X$ with $N=f_{k-1}(X)$
facets.  It was observed by G. Kalai, S. Weinberger and the author that the torsion subgroup  $H_{k-2}(X)_{tor}$
satisfies $|H_{k-2}(X)_{tor}| \leq \sqrt{k}^N$. Kalai on the other hand showed \cite{K83} that there exist $X$'s with $|\thh_{k-2}(X)_{tor}| \geq \sqrt{k/e}^N$. Computer experiments indicate that the $\Rat$-acyclic sum complexes obtained by taking $|A|=k$ often have large torsion. For example, $A=\{0,1,19\} \subset \FF_{83}$ satisfies $|H_1(X_{A,3})| > 1.17^N$ where $N=f_2(X_{A,3})=\binom{82}{2}$. Note that the base of the exponent
$1.17$ is slightly bigger than the constant $\sqrt{3/e}\doteq 1.05$  in Kalai's lower bound.
It view of this it would be interesting to determine (or estimate) the maximum torsion of sum complexes.

\item
Theorem \ref{unhom} characterizes the uncertainty number $u_{\FF}(A)$ with $A \subset K=\FF_p$ and $\FF$ algebraically closed, in terms of the homology of $X_{A,k}$ over $\FF$. It would be useful to find appropriate extensions of this characterization to general finite groups $K$ and arbitrary fields $\FF$.

\end{itemize}

\ \\ \\
{\bf ACKNOWLEDGMENTS}
\\
Research supported by a grant from the Israel Science Foundation
with additional partial support from ERC Advanced Research Grant no 267165 (DISCONV) and
ESF grant (ACAT).
\\
The author would like to thank Shmuel Weinberger for helpful discussions.

\end{document}